\documentclass{article}
\usepackage{fullpage,amssymb}

\newenvironment{@abssec}[1]{%
     \if@twocolumn
       \section*{#1}%
     \else
       \vspace{.05in}\footnotesize
       \parindent .2in
         {\bfseries #1. }\ignorespaces
     \fi}
     {\if@twocolumn\else\par\vspace{.1in}\fi}
\newenvironment{keywords}{\begin{@abssec}{Key words}}{\end{@abssec}}
\newenvironment{AMS}{\begin{@abssec}{AMS subject classification}}{\end{@abssec}}
\newenvironment{proof}{\noindent{\bf Proof.\/}}{\hfill$\Box$}

\newtheorem{prop}{Proposition}
\newtheorem{thm}[prop]{Theorem}
\newtheorem{lem}[prop]{Lemma}
\newtheorem{corol}[prop]{Corollary}

\def\ii{{\rm i}}
\def\st#1{\langle #1 \rangle}
\def\eqbd{\mathop{{:}{=}}}
\def\C{{\rm C\kern-.48em\vrule width.06em height.6em depth-.02em 
                 \kern.48em}}
\def\tr{\mathop{\rm tr}\nolimits}

\begin{document}

\title{$M$-matrices satisfy Newton's inequalities}
\author{ Olga Holtz\thanks{On leave from CS Department, Univ. of Wisconsin, Madison, 
WI 53706, USA. 
Supported by Alexander von Humboldt Foundation.} \\
Institut f\"ur Mathematik, MA 4-5 \\
 Technische Universit\"at Berlin \\
 D-10623 Berlin, Germany \\
holtz@math.TU-Berlin.DE}
\date{}
\maketitle

\begin{keywords}{$M$-matrices, Newton's inequalities, immanantal inequalities, 
generalized matrix functions, quadratic forms, binomial identities, 
nonnegative inverse eigenvalue problem.} \end{keywords}

\begin{AMS}{Primary 15A42; Secondary 15A15, 15A45, 15A48, 15A63. 
05E05, 05A10, 05A17, 05A19, 26D05, 65F18}\end{AMS}

\begin{abstract} Newton's inequalities $c_n^2 \ge c_{n-1}c_{n+1}$ are shown
to hold for the normalized coefficients $c_n$ of the characteristic 
polynomial of any $M$- or inverse $M$-matrix. They are derived by establishing
first an auxiliary set of inequalities also valid for both of these classes. 
They are also used to derive some new necessary conditions on the eigenvalues 
of nonnegative matrices. 
\end{abstract} 

\maketitle

\section{Introduction}

The goal of the paper is to prove a conjecture made in~\cite{HS} about a set of
inequalities satisfied by (the elementary symmetric functions of) the eigenvalues
of any $M$- or inverse $M$-matrix. 

Let $\st{n}$ denote the collection of all increasing sequences with elements
from the set $\{1,2,\ldots,n\}$, let $\# \alpha$ denote the size of the
sequence $\alpha$, and let $\alpha'$ denote the complementary or `dual' 
sequence whose elements are all the integers from $\{1,2,\ldots,n\}$ not in 
$\alpha$. Given a matrix $A\in \C^{n\times n}$, 
the notation $A(\alpha)$ ($A[\alpha]$) will be used for the principal
submatrix (minor) of $A$ whose rows and columns are indexed by $\alpha$.
By convention, $A[\emptyset]\eqbd 1$. 

A matrix $A$ is called a $P$-matrix if $A[\alpha]>0$
for all $\alpha\in \st{n}$. $A$ is called a (nonsingular) $M$-matrix if it 
is a $P$-matrix and its 
off-diagonal entries are nonpositive. If in this definition the positivity 
of all principal minors is relaxed to nonnegativity, one obtains the 
class of all $M$-matrices, including the singular ones. The class of
inverse $M$-matrices consists of matrices whose inverses are $M$-matrices.
The $M$-matrices are an important class arising in many contexts (see, for 
example,~\cite[Chapter 6]{BP}).

Given a matrix $A$, let $c_j(A)$ denote the normalized coefficients
of its characteristic polynomial:
\[ c_j(A)\eqbd\sum_{\# \alpha=j} A[{\alpha}]/{n \choose j}, \qquad
 \qquad  j=0,\ldots,n. \]
The inequalities
\begin{equation} c_j^2(A)\geq c_{j-1}(A)c_{j+1}(A), \qquad \qquad 
j=1,\ldots,n-1,
 \label{newton} \end{equation}
are known for real diagonal matrices, i.e., simply for
sequences of real numbers (see~\cite{Nic} and references therein),  as was
first proved by Newton. Since the numbers $c_j$ are invariant under
similarity, Newton's inequalities~(\ref{newton}) also hold for all 
diagonalizable matrices with real spectrum, and therefore also for the
closure of this set, viz. for {\em all\/} matrices with real spectrum. 

It was conjectured in~\cite{HS} that Newton's inequalities are also satisfied
by $M$- and inverse $M$-matrices (and by matrices similar to those). 
The next section contains proofs of several results on $M$-matrices and 
symmetric functions culminating in the proof of this fact.

\section{Proof of Newton's inequalities}

Let us begin by establishing a set of auxiliary inequalities first.
Given an $n\times n$-matrix $A$ and nonnegative integers $m_1$, $m_2$, $k$, 
define functions $S_{m_1,m_2,k}$ as follows
\begin{equation} S_{m_1,m_2,k}(A)\eqbd \sum_{\alpha \in \st{n}, \# \alpha=m_1, \atop 
\beta\in \st{n}, \# \beta=m_2, \# \alpha\cap \beta=k } A[\alpha] A[\beta]. \label{def}
\end{equation}

\begin{thm} For any  $M$- or inverse $M$-matrix $A$ of order $n$ and 
nonnegative integers $m<n$, $k<m$, 
\begin{equation} S_{m,m,k}(A)/S_{m,m,k}(I_n) \geq S_{m+1,m-1,k}(A)/S_{m+1,m-1,k}(I_n), 
\label{genimm} 
\end{equation} 
where $I_n$ denotes the identity matrix of order $n$.
\end{thm}

\begin{proof} by induction. 

Case 1 (induction base). If $k=0$, $n=2m$, then~(\ref{genimm}) is a special case of 
Theorem~1.3 from~\cite{JJP}. Indeed, since $n=2m$, the functions $S_{m.m,0}$
and $S_{m+1,m-1,0}$ are immanants, $\lambda \eqbd (m,m)$ and $\mu\eqbd (m+1,m-1)$ are 
partitions of $n$, and $\mu$ majorizes $\lambda$. Then the normalized
immanant corresponding to $\mu$ does not exceed the one corresponding to $\lambda$ 
(beware a typo in~\cite{JJP}, where the sign is reversed). If an $M$-matrix $A$ is
nonsingular, then $A^{-1} [\alpha]=A[\alpha']/\det A$ 
(see, e.g.,~\cite[Section 1.4]{G}), hence $S_{m,m,0}(A^{-1})=S_{m,m,0}(A)/(\det A)^2$,
$S_{m+1,m-1,0}(A^{-1})=S_{m+1,m-1,0}(A)/(\det A)^2$, so the inequality~(\ref{genimm})
holds for the matrix $A^{-1}$ as well.

Now assume~(\ref{genimm}) holds for all $M$- and inverse $M$-matrices of order smaller
than $n$. 

Case 2 (induction step of the first kind). Suppose $2m-k<n$ and $A$ is an 
$M$- or inverse $M$-matrix. Then both 
normalized functions $S_{m,m,k}(A)/S_{m,m,k}(I_n)$ and $S_{m+1,m-1,k}(A)/S_{m+1,m-1,k}(I_n)$ 
can be obtained by first averaging the terms $A[\alpha] A[\beta]$ over 
submatrices of order $n-1$ and then taking the average of the obtained $n$ 
quantities:  
\begin{eqnarray*}
{S_{m,m,k}(A)\over S_{m,m,k}(I_n)} &=&{1\over n}\sum_{\alpha\in \st{n}, \# \alpha=n-1} 
{S_{m,m,k}(A(\alpha))\over S_{m,m,k}(I_{n-1})} \\
{S_{m+1,m-1,k}(A) \over S_{m+1,m-1,k}(I_n)}&=&{1\over n}\sum_{\alpha\in \st{n}, \# \alpha=n-1} 
{S_{m+1,m-1,k}(A(\alpha)) \over S_{m+1,m-1,k}(I_{n-1})}. 
\end{eqnarray*}
But principal submatrices of $M$- (inverse $M$-) matrices are again $M$- (inverse $M$-) matrices
(\cite[p.113, p.119]{HJ}), therefore the inductive assumption holds for all submatrices
$A(\alpha)$, $\#\alpha=n-1$. This implies~(\ref{genimm}) for the matrix $A$ itself.

Case 3 (induction step of the second kind). Let $2m-k=n$ and $k>0$. First 
assume $A$ is a {\em nonsingular\/} $M$- or inverse $M$-matrix.
Switch to the dual case: Each $A[\alpha] A[\beta]$ in the right-hand side of~(\ref{def})  
equals $A^{-1}[\alpha']A^{-1}[\beta']/(\det A)^2$, the index
sets $\alpha'$ and $\beta'$ do not intersect, and
$ \#\alpha'+ \#\beta'=2(n-m)<n. $ 
Hence 
\[ S_{m,m,k}(A)={S_{n-m,n-m,0}(A^{-1})\over (\det A)^2}, \qquad
S_{m+1,m-1,k}(A)={S_{n-m+1,n-m-1,0}(A^{-1}) \over (\det A)^2} \]  
and the functions $S_{n-m,n-m,0}(A^{-1})$, $S_{n-m+1,n-m-1,0}(A^{-1})$ are as 
in Case 2 above.
Thus~(\ref{genimm}) holds for the matrix $A^{-1}$ and hence for the matrix 
$A$. So, the induction step of this kind is now proved for nonsingular 
$M$-matrices and their inverses. But the set of all $M$- matrices is the 
closure of the set of nonsingular $M$-matrices  (see, e.g.,~\cite[p.119]{HJ}), 
which justifies the induction step for singular $M$-matrices as well.

With all possible cases considered, the theorem is proved. \end{proof} 

Now, the theorem can be used to replace  Newton's inequalities by a stronger
(but simpler) set of quadratic inequalities in the variables $A[\alpha]$. 

\begin{lem} \label{first} Let $m\in \{1, \ldots, n\}$ be fixed and let $t(m)$ be the 
column vector  \[t(m)\eqbd (t_\alpha)_{\alpha \in \st{n}, \#\alpha =m}.\]
Let $\Psi_m$ denote the Hermitian form
\small
\begin{equation} t(m) \mapsto t(m)^*\Psi_m t(m)\eqbd \sum_{j=0}^m (
m(n-m)  -(m+1)(n-m+1) {m-j\over m-j+1}) \sum_{\# \alpha= \# \beta =m \atop \# \alpha \cap \beta =j}
\overline{t_\alpha} t_\beta. \label{Psi} \end{equation}
\normalsize
If $\Psi_m$ is nonnegative definite, then the $m$th Newton's 
inequality~(\ref{newton}) holds.
\end{lem}

\begin{proof} Expanding both sides of the $m$th Newton's inequality yields
\begin{eqnarray*} c_m^2(A)&=&\sum_{j=0}^m S_{m,m,j}(A) / {n \choose m}^2,  \\
 c_{m-1}(A) c_{m+1}(A)&=&\sum_{j=0}^{m-1} S_{m+1,m-1,j}(A) / {n \choose m+1} {n \choose m-1}, 
 \end{eqnarray*}
So, the $m$th Newton's inequality is equivalent to 
\begin{equation}
m(n-m) \sum_{j=0}^m S_{m,m,j}(A)\geq  (m+1) (n-m+1)
\sum_{j=0}^{m-1} S_{m+1,m-1,j}(A).  \label{equiv} \end{equation} 
On the other hand, straightforward counting gives
\begin{eqnarray*} S_{m,m,j}(I_n) &= &{n \choose j} {n-j \choose m-j} {n-m \choose m-j}, \\
S_{m+1,m-1,j}(I_n)&=&{n \choose j} {n-j \choose m-j-1} {n-m+1 \choose m-j+1}, 
\end{eqnarray*}
hence the inequalities~(\ref{genimm}) are equivalent to
\[ (m-j)S_{m,m,j}(A)\geq (m-j+1)S_{m+1,m-1,j}(A).\] 
Thus, upon replacing each $S_{m+1,m-1,j}$ in the right-hand side of~(\ref{equiv}) 
by ${(m-j) \over (m-j+1)}S_{m,m,j}$, one obtains a set of inequalities stronger than Newton's. 
Precisely, these stronger inequalities assert that
\[ \sum_{j=0}^m (m(n-m)-(m+1)(n-m+1){m-j \over m-j+1})S_{m.m.j} \geq 0, \]
or, recalling the definitions of $S_{m,m,j}$ and of $\Psi_m$,
\[ a(m)^* \Psi_m a(m)\geq 0 \qquad {\rm where} \quad a(m)\eqbd (A[\alpha])_{
\alpha\in\st{n}, \#\alpha=m}. \] 
So, if $\Psi_m$ is nonnegative definite, then the $m$th
Newton's inequality is satisfied.  \end{proof}

Thus, it remains to prove the following.

\begin{lem} \label{last} With the notation of Lemma~\ref{first}, 
$t(m)^*\Psi_m t(m)\geq 0$ for all $t(m)$ and all $m=1, \ldots, n-1$. 
\end{lem}

\begin{proof} Consider first the Hermitian form  
\[ \Phi_m : t(m) \mapsto
t(m)^* \Phi_m t(m) \eqbd \sum_{j=0}^m j\sum_{\# \alpha= \# \beta =m \atop 
\# \alpha \cap \beta =j} \overline{t_\alpha} t_\beta. \]
The representation matrix 
\[ \left( \# \alpha \cap \beta  \right)_{\alpha, \beta} \] 
of this Hermitian form is the Gramian, with respect to the standard inner 
product,  for the system of vectors $ (v_\alpha)_{\alpha}$  where 
\[ v_\alpha(i) \eqbd \left\{ \begin{array}{ll} 1 & {\rm if} \quad i\in \alpha \\ 0 & 
{\rm otherwise,} \end{array} \right. 
\]
hence is nonnegative definite. Moreover, the vector $e$ of all ones 
(of appropriate length) is an eigenvector
of $\Phi_m$. Now consider a form \[ \widetilde\Phi_m : t(m) \mapsto
t(m)^* \widetilde\Phi_m t(m) \eqbd \sum_{j=0}^m (m-j+1)\sum_{\# \alpha= 
\# \beta =m \atop \# \alpha \cap \beta =j} \overline{t_\alpha} t_\beta. \]
Its representation matrix  is obtained by subtracting  $\Phi_m$ from a 
positive multiple of the Hermitian rank-one matrix $e e^*$ (precisely
 $(m+1)e e^*$), therefore all eigenvalues of 
$\widetilde\Phi_m$  are nonpositive except for  the one corresponding 
to the eigenvector $e$, which is strictly positive. Therefore, by~\cite{B}, 
the Hadamard inverse $\widetilde{\Psi}_m$ of the matrix $\widetilde{\Phi}_m$, 
i.e., the matrix 
\[ \left( {1\over m- \# \alpha \cap \beta +1} \right)_{\alpha, \beta} \] 
is nonnegative definite. Finally, $\Psi_m$ is obtained from $(m+1)(n-m+1)
\widetilde{\Psi}_m$ by  subtracting the rank-one matrix $e e^*$ this time 
multiplied by $(n+1)$. The eigenvalue of $\Psi_m$ corresponding to $e$ is 
equal to zero, since
\small
\begin{eqnarray*}
e^*\Psi_m e&=&m(n-m)\sum_{j=0}^m S_{m,m,j}(I_n)-(m+1)(n-m+1)\sum_{j=0}^m{m-j\over m-j+1} 
S_{m,m,j}(I_n) \\
&=& m(n-m)\sum_{j=0}^m S_{m,m,j}(I_n)-(m+1)(n-m+1)\sum_{j=0}^{m-1}S_{m+1,m-1,j}(I_n)=0.
\end{eqnarray*} \normalsize        
All the other eigenvalues of $\Psi_m$ are nonnegative, so $\Psi_m$ is 
nonnegative definite.   
\end{proof}

This lemma finishes the proof of Newton's inequalities.

\begin{thm} \label{main}
Let $A$ be similar to an $M$- or inverse $M$-matrix. Then the normalized 
coefficients of its characteristic polynomial satisfy Newton's 
inequalities~(\ref{newton}).
\end{thm}

Also note that a by-product of Lemma~\ref{last} is a binomial identity:

\begin{corol}
$ \sum_{j=0}^m (m(n-m) - (m+1)(n-m+1) {m-j\over m-j+1}) {m \choose j} 
{n-m \choose m-j}=0. $
\end{corol}

\section{Newton's inequalities and the inverse eigenvalue problem for 
nonnegative matrices}

As possible applications of Theorem~\ref{main} one can envision eigenvalue 
localization for $M$- and inverse $M$-matrices as well as inverse 
eigenvalue problems. In the rest of the paper the focus will be on the 
latter problem for nonnegative matrices. 

The nonnegative inverse eigenvalue problem (NIEP) is that of determining 
necessary and sufficient conditions in order that a given $n$-tuple be the 
spectrum of an entrywise nonnegative $n\times n$ matrix. For details and
history of the problem, see~\cite{BP},~\cite{L},~\cite{M}, and references 
therein. 

Two known necessary conditions that an $n$-tuple 
$\Lambda\eqbd(\lambda_1,\ldots,\lambda_n)$ be realizable
as a spectrum of a nonnegative matrices are formulated in
terms of its moments
$$ s_k(\Lambda)\eqbd \sum_{j=1}^n \lambda_j^k,$$
viz.
\begin{eqnarray}
& s_k\geq 0, & {\rm all} \quad k \label{moments} \\
& s_k^m \leq n^{m-1} s_{km}, & {\rm all} \quad k, m. \label{JLL}
\end{eqnarray}
The condition~(\ref{moments}) follows simply from the fact that
$\tr(A^k)$ is the $k$th moment of the eigenvalue sequence of $A$,
while the condition~(\ref{JLL}) is due to Loewy and London~\cite{LL}
and, independently, Johnson~\cite{J}. 

Newton's inequalities proven above result in a third set of conditions 
necessary for realizability of a given $n$-tuple. Precisely, if
$\Lambda=(\lambda_1,\ldots,\lambda_n)$ is the spectrum of a nonnegative
matrix $A$ and $\lambda_1=\max |\Lambda|$ is its spectral radius, then
the set $(0,\lambda_1-\lambda_2,\ldots,\lambda_1-\lambda_n)$
is the spectrum of an $M$-matrix $\lambda_1 I -A$ and should therefore
satisfy Newton's inequalities~(\ref{newton}).
 
Newton's inequalities are independent of~(\ref{moments}) and~(\ref{JLL}).
First of all, it is clear that~(\ref{newton}) and~(\ref{moments}) are 
independent:  for example, the triple $(1,-1, -1)$ does not 
satisfy~(\ref{moments}) but its shifted counterpart $(0,2,2)$ 
satisfies~(\ref{newton}), while the triple $(\sqrt{2},\ii, -\ii)$
satisfies~(\ref{moments}) but the corresponding shifted triple 
$(0,\sqrt{2}-i, \sqrt{2}+i)$ does not satisfy~(\ref{newton}).

Moreover, neither the two conditions~(\ref{moments}) and~(\ref{newton})
{\it together\/} imply~(\ref{JLL}) nor the two conditions~(\ref{moments}) 
and~(\ref{JLL}) {\it together\/} imply~(\ref{newton}). 

Indeed, the conditions~(\ref{moments}) 
and~(\ref{newton}) can be satisfied while the conditions~(\ref{JLL})
may fail. To show this, consider the 10-tuple 
$\Lambda\eqbd(3,1,1,1,1,1,-2,-2,-2,-2).$
Its first and third moment are equal to zero, while the rest are positive.
Now, let us introduce its perturbed version 
$\Lambda_t\eqbd (3+t_1,1+t_2,1,1,1,1,-2+t_3,-2,-2,-2)$, where
the $t$'s are real and 
\begin{eqnarray*}
&& t_1+t_2+t_3>0 \\
&& (3+t_1)^3+(1+t_2)^3+(-2+t_3)^3=20, 
\end{eqnarray*}
which is always possible according to the Linearization Lemma~\cite[p.163]{Ma},
since the system
\begin{eqnarray*}
&& t_1+t_2+t_3>0 \\
&& 9t_1+t_2+4t_3=0 
\end{eqnarray*}
is solvable arbitrarily close to the point $(0,0,0)$. 
The first moment of $\Lambda_t$ is thus positive, while the
third is still zero. All the other moments remain positive if $(t_1,t_2,t_3)$
is sufficiently small. So, ~(\ref{moments}) is satisfied. The Newton
conditions~(\ref{newton}) are satisfied as well, since $\Lambda_t$ is real.
But the condition~(\ref{JLL}) with $k=1$, $m=3$ fails. 

To construct an example where~(\ref{moments}) and~(\ref{JLL}) are satisfied 
but~(\ref{newton}) fails, consider the sequence of zeros of the polynomial 
$p(x)=x^6 -6x^5 + 14x^4-20x^3$. It does not satisfy~(\ref{newton}): This
polynomial is obtained by cutting the expansion of $(x-1)^6$, 
whose coefficients satisfy~(\ref{newton}) with strict equalities, 
and then decreasing slightly the value (originally $15$) of one coefficient.
Then the second Newton's inequality fails. The non-zero roots of $p$
are approximately $3.6702$ and  $1.1649 \pm  2.0229\ii$. By shifting back 
by the largest absolute value $\approx 3.6702$, one obtains the $6$-tuple 
$\Lambda \eqbd(a,a,a,0, b, \overline{b})$ with $a\approx 3.6702$, 
$b\approx  2.5054 + 2.0229\ii$.
It is not hard, though a bit tedious, to check that $\Lambda$ 
satisfies~(\ref{JLL}). Since $s_1(\Lambda)>0$, this also implies 
that all moments of $\Lambda$ are positive. This shows that~(\ref{newton}) 
cannot be derived from~(\ref{moments}) and~(\ref{JLL}).

In the case the first moment of an $n$-tuple is zero, Laffey and 
Meehan~\cite{LM} established another necessary condition, viz.,
$$ (n-1) s_4 \geq s_2^2.  $$
It is also not implied by~(\ref{moments}),~(\ref{JLL}) and~(\ref{newton}).
An example is provided by the $5$-tuple $(3,3,-2,-2,-2)$.

Note, however, that the condition~(\ref{JLL}) 
with $k=1$, $m=2$ is exactly equivalent to the Newton 
inequality~(\ref{newton}) for $j=1$.

\section*{Acknowledgements} I am grateful to Hans Schneider, Thomas Laffey,
and an anonymous referee for helpful remarks and suggestions.

\end{document}